# A COMPLEMENTARY DESIGN THEORY FOR DOUBLING

By Hongquan Xu[1] and Ching-Shui Cheng[2]

*University of California, Los Angeles, and University of California, Berkeley*

Chen and Cheng [*Ann. Statist.* **34** (2006) 546–558] discussed the method of doubling for constructing two-level fractional factorial designs. They showed that for $9N/32 \leq n \leq 5N/16$, all minimum aberration designs with $N$ runs and $n$ factors are projections of the maximal design with $5N/16$ factors which is constructed by repeatedly doubling the $2^{5-1}$ design defined by $I = ABCDE$. This paper develops a general complementary design theory for doubling. For any design obtained by repeated doubling, general identities are established to link the wordlength patterns of each pair of complementary projection designs. A rule is developed for choosing minimum aberration projection designs from the maximal design with $5N/16$ factors. It is further shown that for $17N/64 \leq n \leq 5N/16$, all minimum aberration designs with $N$ runs and $n$ factors are projections of the maximal design with $N$ runs and $5N/16$ factors.

**1. Introduction.** Fractional factorial designs are widely used in various experiments and are typically chosen according to the *minimum aberration* criterion. There is much recent work on the theory and construction of minimum aberration designs; see Wu and Hamada [11] for a good summary. Chen and Cheng [4] discussed the method of doubling for constructing two-level fractional factorial designs, in particular, those of resolution IV. Suppose that $X$ is an $N \times n$ matrix with two distinct entries, 1 and $-1$. Then the *double* of $X$ is the $2N \times 2n$ matrix $\begin{pmatrix} X & X \\ X & -X \end{pmatrix}$. Doubling plays an important role in the construction of maximal designs. A regular design of resolution IV or higher is called *maximal* if its resolution reduces to three whenever an extra factor is added. It is obvious that every regular design of resolution

Received August 2006; revised February 2007.
[1]Supported in part by NSF Grant DMS-05-05728.
[2]Supported in part by NSF Grant DMS-05-05556.

IV or higher is a projection of some maximal regular design of resolution IV or higher. Some recent results in the literature of finite projective geometry essentially characterize all maximal regular designs of resolution IV with $n \geq N/4 + 1$. These results imply that for $n \geq N/4 + 1$, a maximal regular design of resolution IV or higher must have

$$n \in \{N/2, 5N/16, 9N/32, 17N/64, 33N/128, \ldots\},$$

where $N$ is a power of 2 and $n$ is the number of factors; see Chen and Cheng [4]. All these maximal regular designs can be constructed by repeatedly doubling some small designs. For example, the maximal regular design with $N/2$ factors can be constructed by repeatedly doubling the $2^1$ design. The maximal regular design with $5N/16$ factors is constructed by repeatedly doubling the $2^{5-1}$ design defined by $I = ABCDE$. An immediate result is that all regular designs of resolution IV with $5N/16 < n \leq N/2$ must be projections of the maximal regular design with $N/2$ factors and that all regular designs of resolution IV with $9N/32 < n \leq 5N/16$ must be projections of the maximal regular design with either $N/2$ or $5N/16$ factors. Chen and Cheng [4] showed that for $9N/32 \leq n \leq 5N/16$, all minimum aberration designs with $N$ runs and $n$ factors are projections of the maximal design with $5N/16$ factors. The motivation of this paper is to develop a complementary design theory that guides the construction of minimum aberration designs from the maximal design with $5N/16$ factors. Previously, Chen and Hedayat [6] and Tang and Wu [10] developed a complementary design theory for saturated regular designs of resolution III. Butler [3] and Chen and Cheng [5] developed a complementary design theory for the maximal regular designs with $N/2$ factors. Note that the maximal regular designs with $N/2$ factors are even designs (i.e., all defining words have even lengths) and are the only even designs that are maximal. We shall call them maximal even designs.

In Section 2, we review some basic concepts and background. In Section 3, we develop a general complementary design theory for doubling. For each design constructed by repeated doubling, we express the wordlength pattern of a projection design in terms of a linear combination of the wordlength pattern of its complement and some terms that depend on the frequencies of the columns being doubled. We illustrate how previous complementary design theories for saturated regular designs of resolution III and maximal even designs can be easily derived as special cases of our result. In Section 4, we apply the general theory to the maximal design with $5N/16$ factors and obtain a general rule on constructing its minimum aberration projection designs. As conjectured by Chen and Cheng [4], we show that the minimum aberration projection designs have equal or nearly equal frequencies of the columns being doubled. In Section 5, we consider the maximal design with $9N/32$ factors and derive a lower bound for the wordlength patterns of its projections. This lower bound is then used in Section 6 to show that for $17N/64 \leq n \leq 5N/16$, all minimum aberration designs with $N$ runs and $n$ factors are projections of the maximal design with $5N/16$ factors.

**2. Basic concepts and background.** An $N$-run design with $n$ factors is represented by an $N \times n$ matrix, where each row corresponds to a run and each column corresponds to a factor. A two-level design takes on only two symbols which, for convenience, are denoted by 0 and 1 in the rest of the paper. Regular designs are those which can be constructed by using defining relations and are discussed in many textbooks on experimental design; see, for example, [9, 11]. A regular $2^{n-p}$ design is defined by $p$ independent defining words and has $N = 2^{n-p}$ runs and $n$ factors. The $p$ independent defining words together generate $2^p - 1$ defining words. The *resolution* is the length of the shortest defining word. Let $A_i$ be the number of defining words of length $i$. Then the vector $(A_1, \ldots, A_n)$ is called the *wordlength pattern*. The *minimum aberration* criterion introduced by Fries and Hunter [7] chooses a design by sequentially minimizing $A_1, A_2, \ldots, A_n$. Throughout the paper, we will consider only two-level regular designs.

For a two-level regular $N \times n$ design $D$ (with symbols 0 and 1), let $w_i(D)$ be the *Hamming weight* (i.e., the number of nonzero components) of the $i$th row of $D$. For integers $k > 0$, define the moments $M_k(D) = \sum_{i=1}^{N} w_i(D)^k$. Using the Pless power moment identities [8], a fundamental result in coding theory, Xu [12, 13, 14] expressed $M_k$ as a linear combination of $A_1, \ldots, A_k$ as follows.

LEMMA 1. *For a two-level regular $N \times n$ design $D$ and integers $k > 0$,*

$$M_k(D) = N \sum_{i=0}^{\min(n,k)} Q_k(i;n) A_i(D),$$

*where* $Q_k(i;n) = (-1)^i \sum_{j=0}^{k} j! S(k,j) 2^{-j} \binom{n-i}{j-i}$, $S(k,j) = (1/j!) \times \sum_{i=0}^{j} (-1)^{j-i} \binom{j}{i} i^k$ *is a Stirling number of the second kind and* $A_0(D) = 1$.

It is easy to verify that $S(k,k) = 1$ and $Q_k(k;n) = (-1)^k k!/2^k$.

**3. Doubling and complementary designs.** For a two-level regular $N_0 \times m_0$ design $X_0 = (x_{ij})$ with $x_{ij} = 0$ or 1. Define its double as

$$\begin{pmatrix} X_0 & X_0 \\ X_0 & X_0 + 1 \end{pmatrix} \pmod{2}.$$

Let $X$ be the regular $N \times m$ design which is obtained by repeatedly doubling $X_0$ $t$ times, where $N = 2^t N_0$ and $m = 2^t m_0$.

Suppose that $D$ and $\overline{D}$ are a pair of complementary projection designs of $X$, that is, $D \cup \overline{D} = X$, where $D$ has $n$ columns and $\overline{D}$ has $u$ columns, with $n + u = m$. Note that each column of $X_0$ generates $2^t$ columns of $X$. For $i = 1, \ldots, m_0$, let $f_i$ be the number of columns of $\overline{D}$ that are generated from the $i$th column of $X_0$. Clearly $\sum_{i=1}^{m_0} f_i = u$.

It is easy to verify the relationship of the Hamming weights of $D$ and $\overline{D}$ as follows:

$$w_i(\overline{D}) = \sum_{j=1}^{m_0} f_j x_{ij} \qquad \text{for } i = 1, \ldots, N_0,$$

$$w_i(D) = \sum_{j=1}^{m_0} (2^t - f_j) x_{ij} = 2^t w_i(X_0) - w_i(\overline{D}) \qquad \text{for } i = 1, \ldots, N_0,$$

$$w_i(D) + w_i(\overline{D}) = m/2 \qquad \text{for } i = N_0 + 1, \ldots, N.$$

Then

$$M_k(D) = \sum_{i=1}^{N} w_i(D)^k = \sum_{i=1}^{N} [m/2 - w_i(\overline{D})]^k + \Delta_k(D, \overline{D})$$

$$= \sum_{i=1}^{N} \sum_{j=0}^{k} \binom{k}{j} (m/2)^{k-j} (-1)^j w_i(\overline{D})^j + \Delta_k(D, \overline{D})$$

$$= \sum_{j=0}^{k} \binom{k}{j} (m/2)^{k-j} (-1)^j M_j(\overline{D}) + \Delta_k(D, \overline{D}),$$

where

(1) $$\Delta_k(D, \overline{D}) = \sum_{i=1}^{N_0} [w_i(D)^k - (m/2 - w_i(\overline{D}))^k]$$

depends on the frequencies $f_1, \ldots, f_{m_0}$ and $X_0$ only. Using Lemma 1, we obtain the following relationship between the wordlength patterns of $D$ and $\overline{D}$.

THEOREM 1. *Let $D$ and $\overline{D}$ be a pair of complementary projection designs of the regular $N \times m$ design $X$ constructed via doubling. For integers $k > 0$,*

$$\sum_{i=0}^{k} Q_k(i; n) A_i(D) = \sum_{i=0}^{k} \left[ \sum_{j=i}^{k} \binom{k}{j} \left(\frac{m}{2}\right)^{k-j} (-1)^j Q_j(i; u) \right] A_i(\overline{D}) + \frac{\Delta_k(D, \overline{D})}{N},$$

*where $Q_k(i; n)$ is defined in Lemma 1, $\Delta_k(D, \overline{D})$ is defined in (1), $A_0(D) = A_0(\overline{D}) = 1$, $A_i(D) = 0$ for $i > n$ and $A_i(\overline{D}) = 0$ for $i > u$.*

Because $Q_k(k;n) = (-1)^k k!/2^k$, we can express $A_k(D)$ in terms of a linear combination of $A_k(\overline{D}), \ldots, A_1(\overline{D})$, plus some terms that depend on the frequencies $f_1, \ldots, f_{m_0}$ and $X_0$ as follows:

$$A_k(D) = (-1)^k A_k(\overline{D}) + c_{k-1} A_{k-1}(\overline{D}) + \cdots + c_1 A_1(\overline{D}) + c_0$$
(2)
$$+ d_k \Delta_k(D, \overline{D}) + \cdots + d_1 \Delta_1(D, \overline{D}),$$

where $c_i$ and $d_i$ are constants that do not depend on $D$ or $\overline{D}$. If $X$ has resolution $r$, then $A_i(D) = A_i(\overline{D}) = 0$ for $1 \leq i < r$ and Theorem 1 implies that $\Delta_i(D, \overline{D})$ is a constant for $1 \leq i < r$. In this paper, we are interested in resolution IV designs. The next result follows from Theorem 1.

COROLLARY 1. *If $X$ has resolution IV or higher, then $\Delta_1(D, \overline{D}) = 0$, $\Delta_2(D, \overline{D}) = N(2n-m)/4$, $\Delta_3(D, \overline{D}) = 3Nn(2n-m)/8$ and*

$$A_4(D) = A_4(\overline{D}) - (2n-m)(6n^2 + 3m - 2)/24 + 2\Delta_4(D, \overline{D})/(3N).$$

The following two examples illustrate how complementary design theories for saturated regular designs of resolution III and maximal even designs can be derived from Theorem 1.

EXAMPLE 1. Let

$$X_0 = \begin{pmatrix} 0 & 0 \\ 0 & 1 \end{pmatrix}.$$

Repeatedly doubling $X_0$ $t$ times yields an $N \times m$ design $X$, where $N = m = 2^{t+1}$. Note that the first column of $X$ is a vector of 0's. Deleting the first column, we obtain a saturated regular design of resolution III, denoted by $X'$. Following the previous notation, we have $w_1(\overline{D}) = w_1(D) = 0$, $w_2(\overline{D}) = f_2$ and $w_2(D) = 2^t - f_2$. Then

$$\Delta_k(D, \overline{D}) = \sum_{i=1}^{2} [w_i(D)^k - (m/2 - w_i(\overline{D}))^k]$$
$$= (2^t - f_2)^k - (2^t)^k - (2^t - f_2)^k = -2^{tk},$$

which does not depend on the frequencies $f_1$ and $f_2$ at all. By (2), sequentially minimizing $A_1(D), A_2(D), \ldots$ is equivalent to sequentially minimizing $-A_1(\overline{D}), A_2(\overline{D}), -A_3(\overline{D}), A_4(\overline{D}), \ldots$. Note that $A_1(\overline{D})$ is maximized if and only if $\overline{D}$ contains the first column of $X$, a vector of 0's. Let $\overline{D}'$ be the complement of $D$ with respect to $X'$, the saturated design of resolution III. It is easy to see that $A_1(\overline{D}') = A_2(\overline{D}') = 0$ and $A_i(\overline{D}) = A_i(\overline{D}') + A_{i-1}(\overline{D}')$ for $1 \leq i \leq u$. Therefore, sequentially minimizing $A_1(D), A_2(D), \ldots$ is equivalent to sequentially minimizing $-A_3(\overline{D}'), A_4(\overline{D}'), -A_4(\overline{D}'), A_5(\overline{D}'), \ldots$. This result is equivalent to the rules developed by Chen and Hedayat [6] and Tang and Wu [10].

EXAMPLE 2. Let $X_0 = \binom{0}{1}$. Repeatedly doubling $X_0$ $t$ times yields an $N \times m$ design $X$, where $N = 2^{t+1}$ and $m = 2^t$. According to Chen and Cheng [4], $X$ is a maximal even design. Following the previous notation, we have $w_1(\overline{D}) = w_1(D) = 0$, $w_2(\overline{D}) = u$ and $w_2(D) = 2^t - u$. Then

$$\Delta_k(D, \overline{D}) = \sum_{i=1}^{2}[w_i(D)^k - (m/2 - w_i(\overline{D}))^k] = (2^t - u)^k - (2^{t-1})^k - (2^{t-1} - u)^k,$$

which is a constant. Since $X$ is an even design, $A_i(D) = A_i(\overline{D}) = 0$ for all odd $i$. Then by (2), sequentially minimizing $A_4(D), A_6(D), \ldots$ is equivalent to sequentially minimizing $A_4(\overline{D}), A_6(\overline{D}), \ldots$. Combining this with the fact that all regular designs of resolution IV with $5N/16 < n \le N/2$ must be projections of the maximal even design $X$, that is, they must be even designs, we reach the conclusion that for $5N/16 < n \le N/2$, $D$ has minimum aberration among all regular designs with $N$ runs and $n$ factors if and only if $\overline{D}$ has minimum aberration among all even designs. This is indeed Theorem 3 of Butler [3]. Chen and Cheng [5] derived explicit relationships between the wordlength patterns of complementary projection designs of a maximal even design.

**4. Complementary projection designs of the maximal regular design of resolution IV with $5N/16$ factors.** Let $X_0$ be the unique $16 \times 5$ design of resolution V defined by $I = ABCDE$. Repeatedly doubling $X_0$ $t$ times, we obtain the maximal regular design of resolution IV with $5N/16$ factors, denoted by $X$, where $N = 16 \cdot 2^t$. Following the previous notation, the Hamming weights of $\overline{D}$ [i.e., $w_i(\overline{D})$] are 0, $f_i + f_j$ $(i < j)$ and $u - f_i$ $(i = 1, \ldots, 5)$. The Hamming weights of $D$ [i.e., $w_i(D)$] are 0, $2^{t+1} - f_i - f_j$ $(i < j)$ and $2^{t+2} - u + f_i$ $(i = 1, \ldots, 5)$. Then

$$\Delta_k(D, \overline{D}) = \sum_{i=1}^{5}[(2^{t+2} - u + f_i)^k - (m/2 - u + f_i)^k]$$
$$+ \sum_{1 \le i < j \le 5}[(2^{t+1} - f_i - f_j)^k - (m/2 - f_i - f_j)^k] - (m/2)^k,$$

where $m = 5 \cdot 2^t$. According to Corollary 1, $\Delta_1(D, \overline{D}) = 0$, $\Delta_2(D, \overline{D}) = 2^{t+2}(2n - 5 \cdot 2^t)$ and $\Delta_3(D, \overline{D}) = 3 \cdot 2^{t+1} n(2n - 5 \cdot 2^t)$, which can be verified easily. For $k \ge 4$, $\Delta_k(D, \overline{D})$ depends on the frequencies $f_1, \ldots, f_5$. The next lemma states that when $u \le 15 \cdot 2^{t-3}$, $\Delta_4(D, \overline{D})$ achieves the minimum value if and only if all $f_i$ are equal or nearly equal.

LEMMA 2. *When $u \le 15 \cdot 2^{t-3}$, $\Delta_4(D, \overline{D})$ achieves the minimum value if and only if pairwise differences of the frequencies $f_1, \ldots, f_5$ do not exceed 1, that is, $|f_i - f_j| \le 1$ for all $1 \le i < j \le 5$.*

PROOF. We show that $\Delta_4(D,\overline{D})$ does not attain the minimum value if $|f_i - f_j| > 1$ for some $i < j$. Without loss of generality, assume that $f_1 \geq \cdots \geq f_5 \geq 0$ and $f_1 - f_5 > 1$. Consider another pair of complementary projection designs $D'$ and $\overline{D}'$ with frequencies $f_1 - 1, f_2, f_3, f_4, f_5 + 1$. It can be verified that

$$\Delta_4(D,\overline{D}) - \Delta_4(D',\overline{D}') = 3 \cdot 2^{t+2}(f_1 - f_5 - 1)(3 \cdot 2^t + 2f_1 + 2f_5 - 2u).$$

Because $f_1 \geq \cdots \geq f_5 \geq 0$ and $u \leq 15 \cdot 2^{t-3}$, we have $f_1 + f_5 > u/5$, $3 \cdot 2^t + 2f_1 + 2f_5 - 2u > 3 \cdot 2^t - 8u/5 \geq 0$ and $\Delta_4(D,\overline{D}) > \Delta_4(D',\overline{D}')$. This completes the proof. □

When $u > 15 \cdot 2^{t-3}$, Lemma 2 does not hold. For example, when $t = 2$ and $u = 8$, $\Delta_k(D,\overline{D})$ has the same value for (a) $f_1 = f_2 = f_3 = 2$ and $f_4 = f_5 = 1$ and (b) $f_1 = f_2 = f_3 = f_4 = 2$ and $f_5 = 0$.

Note that $\Delta_k(D,\overline{D})$ is uniquely determined when the frequencies $f_1, \ldots, f_5$ differ from one another by at most 1. Theorem 1 and Lemma 2 together lead to the following result.

THEOREM 2. *Let $D$ and $\overline{D}$ be a pair of complementary projection designs of the maximal regular design $X$ of resolution IV with $N = 16 \cdot 2^t$ runs and $5N/16$ factors. When $25N/128 \leq n \leq 5N/16$, $D$ has minimum aberration among all possible n-factor projections of $X$ if the following two conditions hold:*

(i) $A_4(\overline{D}), -A_5(\overline{D}), A_6(\overline{D}), -A_7(\overline{D}), \ldots$ *are sequentially minimized.*
(ii) $|f_i - f_j| \leq 1$ *for all $1 \leq i < j \leq 5$.*

The construction of $\overline{D}$ satisfying both conditions (i) and (ii) is challenging. It is possible that there may not exist designs satisfying both conditions. In such cases, Theorem 1 has to be used instead. Here we present complementary designs $\overline{D}$ satisfying both conditions for $u = 1$–11. For $i = 1, \ldots, 5$, the $i$th column of $X_0$ generates columns $5k + i$, $0 \leq k \leq 2^t - 1$, of $X$. For $u = 1, 2, 3, 4$, we can choose $\overline{D}$ as the first $u$ columns of $X$, which satisfies both conditions. For $u = 5$, $\overline{D} = \{1, 2, 3, 4, 5\}$ (i.e., the first five columns of $X$) satisfies both conditions because $A_4(\overline{D}) = 0$ is minimized and $A_5(\overline{D}) = 1$ is maximized. For $u = 6$–11 and $N \geq 128$, let $S = \{1, 2, 3, 4, 5, 6, 12, 18, 24, 30, 31\}$ and $\overline{D}$ be the projection design of $X$ whose column indexes are the first $u$ elements of $S$. Clearly $\overline{D}$ satisfies the second condition. It can be verified (via complete computer search) that $A_4(\overline{D}) = 0$ is minimized and $A_5(\overline{D})$ is maximized among all possible regular designs. Furthermore, except for $u = 9$, such designs are unique up to isomorphism and therefore satisfy both conditions. For $u = 9$ and $N \geq 128$, there are two nonisomorphic designs with $A_4(\overline{D}) = 0$ and maximum $A_5(\overline{D}) = 2$ as follows:

$$\overline{D}_1 = \{1, 2, 3, 4, 5, 6, 12, 18, 24\}, \qquad (A_4, \ldots, A_9) = (0, 2, 1, 0, 0, 0);$$
$$\overline{D}_2 = \{1, 2, 3, 4, 5, 6, 12, 23, 39\}, \qquad (A_4, \ldots, A_9) = (0, 2, 0, 0, 1, 0).$$

The second design has a smaller $A_6$; therefore, we shall choose $\overline{D} = \overline{D}_2$.

**5. Complementary projection designs of the maximal regular design of resolution IV with $9N/32$ factors.** There is a unique maximal regular design of resolution IV with $9N/32$ factors [4], which can be constructed by repeatedly doubling a regular $32 \times 9$ design. For convenience, we choose $X_0$ as the $2^{9-5}$ design with defining contrast subgroup

$$I = 1235 = 2346 = 3457 = 4561 = 5672 = 6713 = 7124 = 123456789$$
$$= 13489 = 24589 = 35689 = 46789 = 57189 = 61289 = 72389.$$

Note that factors 1–7 are cyclic (i.e., replacing factors 1–7 with 2–7, 1 yields the same design).

The explicit expression of $\Delta_k(D, \overline{D})$ is quite complicated. Nevertheless, it is cyclic in $f_1, \ldots, f_7$ and symmetric in $f_8$ and $f_9$. In particular, with some tedious algebra, we have

$$\begin{aligned}\Delta_4(D, \overline{D}) = {}& 9534T^4 - 8T^3(535F_1 + 511G_1) \\ & + 12T^2(51F_1^2 + 3F_2 + 94F_1G_1 + 47G_1^2 + 7G_2) \\ & - 8T(4F_1^3 + 9F_1^2G_1 + 3F_2G_1 \\ & \qquad + 12G_1G_2 + 9F_1G_1^2 + 3F_1G_2 - 8G_3) \\ & + 48T(f_1f_3f_4 + f_2f_4f_5 + f_3f_5f_6 \\ & \qquad + f_4f_6f_7 + f_5f_7f_1 + f_6f_1f_2 + f_7f_2f_3),\end{aligned}$$

where $T = 2^t$, $F_1 = \sum_{i=1}^{7} f_i$, $F_2 = \sum_{i=1}^{7} f_i^2$, $G_1 = f_8 + f_9$, $G_2 = f_8^2 + f_9^2$ and $G_3 = f_8^3 + f_9^3$.

The next lemma gives necessary conditions for $\Delta_4(D, \overline{D})$ to achieve the minimum value and a lower bound.

LEMMA 3. *Let $D$ and $\overline{D}$ be a pair of complementary projection designs of the maximal regular design with $32 \cdot 2^t$ runs and $9 \cdot 2^t$ factors defined in this section. When $u \leq 3 \cdot 2^{t-1}$,*

(a) *necessary conditions for $\Delta_4(D, \overline{D})$ to achieve the minimum value are* (i) $f_8 = f_9 = 0$ *and* (ii) $|f_i - f_j| \leq 1$ *for all $1 \leq i < j \leq 7$;*
(b) $\Delta_4(D, \overline{D}) \geq 2^{t+1}(760n^3 - 5400n^2 2^t + 17380n 2^{2t} - 39477 \cdot 2^{3t})/49.$

PROOF. (a) Let $T = 2^t$ and $F(f_1, \ldots, f_9) = \Delta_4(D, \overline{D})$. Without loss of generality, assume that $f_1$ is the smallest among all $f_i$, $1 \leq i \leq 7$. It can be verified that

$$\begin{aligned}
F(f_1, \ldots, f_9) &- F(f_1 + f_8, f_2, f_3, f_4, f_5, f_6, f_7, 0, f_9) \\
&= 192T^3 f_8 - 24T^2(7f_1 + 4f_2 + 4f_3 + 4f_4 + 4f_5 + 4f_6 + 4f_7)f_8 \\
&\quad + 48T f_8(f_1 f_2 + f_1 f_3 + f_2 f_3 + f_1 f_4 + f_2 f_4 + f_1 f_5 \\
&\qquad + f_2 f_5 + f_3 f_5 + f_4 f_5 + f_1 f_6 + f_3 f_6 + f_4 f_6 \\
&\qquad + f_5 f_6 + f_1 f_7 + f_2 f_7 + f_3 f_7 + f_4 f_7 + f_6 f_7 + f_1 f_9) \\
&> 192T^3 f_8 - 24T^2(5u)f_8 \geq 0
\end{aligned}$$

when $u \leq 3 \cdot 2^{t-1}$. Therefore, a necessary condition for $\Delta_4(D, \overline{D})$ to achieve the minimum value is $f_8 = 0$. By the symmetry of $f_8$ and $f_9$, another necessary condition is $f_9 = 0$. Next, we show that $F(f_1, \ldots, f_9)$ cannot be the minimum if $|f_i - f_j| \geq 2$ for some $1 \leq i < j \leq 7$. Without loss of generality, assume that $f_2 - f_1$ is the largest among all possible $f_i - f_j$ for $1 \leq i, j \leq 7$; for the other cases, the proof is the same due to the special structure of $X_0$. It is sufficient to show that if $f_2 - f_1 \geq 2x > 0$, then

(3) $\qquad F(f_1, \ldots, f_9) > F(f_1 + x, f_2 - x, f_3, f_4, f_5, f_6, f_7, f_8, f_9).$

It is straightforward to verify that

(4)
$$\begin{aligned}
F(f_1, \ldots, f_9) &- F(f_1 + x, f_2 - x, f_3, f_4, f_5, f_6, f_7, f_8, f_9) \\
&= 24Tx[(f_2 - f_1 - x)(3T - 2(f_6 + f_8 + f_9)) - 2(f_3 - f_5)(f_4 - f_7)].
\end{aligned}$$

By the assumption that $f_2 - f_1$ is the largest among all $f_i - f_j$ for $1 \leq i, j \leq 7$, it is obvious that

$$-2(f_3 - f_5)(f_4 - f_7) \geq -(f_2 - f_1)(f_3 + f_4 + f_5 + f_7).$$

If $f_2 - f_1 \geq 2x > 0$, then $f_2 - f_1 - x \geq (f_2 - f_1)/2$ and

(5)
$$\begin{aligned}
(f_2 - f_1 - x)&(3T - 2(f_6 + f_8 + f_9)) - 2(f_3 - f_5)(f_4 - f_7) \\
&\geq ((f_2 - f_1)/2)(3T - 2(f_6 + f_8 + f_9)) \\
&\quad - (f_2 - f_1)(f_3 + f_4 + f_5 + f_7) \\
&= ((f_2 - f_1)/2)(3T - 2(f_6 + f_8 + f_9 + f_3 + f_4 + f_5 + f_7)) \\
&> (f_2 - f_1)(3T - 2u)/2 \geq 0,
\end{aligned}$$

because $u \leq 3T/2 = 3 \cdot 2^{t-1}$. Then (3) follows from (4) and (5). This completes the proof.

(b) Following the proof of (a), when $f_i = u/7 = (9 \cdot 2^t - n)/7$ for $i = 1, \ldots, 7$ and $f_8 = f_9 = 0$, $F(f_1, \ldots, f_9)$ achieves the lower bound. The lower bound is achievable only if $u$ is a multiple of 7. $\square$

A result similar to Theorem 2 can be developed. However, such a result is not of interest because, as we will show in the next section, projection designs of the maximal regular design with $9N/32$ factors often do not have minimum aberration. The following lower bound of $A_4$, a direct result of Corollary 1 and Lemma 3(b), is important in the development.

LEMMA 4. *Let $D$ be an $n$-factor projection design of the maximal regular design of resolution IV with $N = 32 \cdot 2^t$ runs and $9N/32$ factors. For $15N/64 \le n \le 9N/32$,*
$$A_4(D) \ge L(n) = [a(n) - b(n)2^t + c(n)2^{2t} - 39477 \cdot 2^{3t}]/1176,$$
*where $a(n) = 196n + 172n^3$, $b(n) = 882 + 2646n + 2754n^2$ and $c(n) = 11907 + 17380n$.*

**6. Some results on minimum aberration designs.** Chen and Cheng [4] showed that for $9N/32 \le n \le 5N/16$, minimum aberration designs are projection designs of the maximal regular design of resolution IV with $5N/16$ factors. Here we strengthen their result and show that this is also true for $17N/64 \le n < 9N/32$.

We first establish a general upper bound of wordlengths for projection designs.

LEMMA 5. *Let $X$ be a regular design of resolution $r$ with $m$ factors and $A_r(X)$ words of length $r$. For $r \le n \le m$, there exists an $n$-factor projection design having at most $A_r(X)\binom{n}{r}/\binom{m}{r}$ words of length $r$.*

PROOF. We show that such a design can be constructed by deleting one factor at a time. First, there must be a factor that appears in at least $r \cdot A_r(X)/m$ words of length $r$. Deleting this factor yields a design with $m - 1$ factors and at most $A_r(X) - r \cdot A_r(X)/m = A_r(X)(m-r)/m$ words of length $r$. Repeat this procedure until we get an $n$-factor projection design which has at most
$$A_r(X)\frac{m-r}{m}\frac{m-1-r}{m-1}\cdots\frac{n+1-r}{n+1} = A_r(X)\binom{n}{r}\bigg/\binom{m}{r}$$
words of length $r$. $\square$

For the maximal regular design $X$ of resolution IV with $5N/16 = 5 \cdot 2^t$ factors, Chen and Cheng [4] showed that $A_4(X) = (65 \cdot 2^{3t-2} - 75 \cdot 2^{2t-2} + 5 \cdot 2^{t-1})/6$. Applying Lemma 5, we obtain the following upper bound of $A_4$ for minimum aberration designs.

COROLLARY 2. *For $N = 16 \cdot 2^t$, $t \geq 0$, and $n \leq 5N/16$, there exists a regular $N \times n$ design $D$ of resolution IV with*

$$A_4(D) \leq ((65 \cdot 2^{3t-2} - 75 \cdot 2^{2t-2} + 5 \cdot 2^{t-1})/6) \binom{n}{4} \bigg/ \binom{5 \cdot 2^t}{4}.$$

Corollary 2 is used next to show that for $N/4 + 1 \leq n \leq 5N/16$, projection designs of the maximal regular designs of resolution IV with $N/2$ or $9N/32$ factors do not have minimum aberration.

LEMMA 6. *For $N = 16 \cdot 2^t$, $t \geq 0$, and $N/4 + 1 \leq n \leq 5N/16$, projection designs of the maximal regular design of resolution IV with $N/2$ factors do not have minimum aberration.*

PROOF. Let $D$ be an $n$-factor projection design of the maximal regular design of resolution IV with $N/2 = 2^{t+3}$ factors. It follows from (2.2) and (4.3) of [4] that

$$A_4(D) \geq \left[ \binom{n}{2}^2 \bigg/ (2^{t+3} - 1) - \binom{n}{2} \right] \bigg/ 6.$$

It is sufficient to show that the upper bound given in Corollary 2 is less than the above lower bound. It is equivalent to show that

$$[(65 \cdot 2^{3t-2} - 75 \cdot 2^{2t-2} + 5 \cdot 2^{t-1})](n-2)(n-3)$$
$$(6) \qquad < \left[ \frac{6n(n-1)}{2^{t+3} - 1} - 12 \right] \binom{5 \cdot 2^t}{4}.$$

Given $t$, both the left-hand side and the right-hand side are polynomials in $n$ of degree 2. It can easily be verified that (6) holds for $t \geq 0$ and $4 \cdot 2^t + 1 \leq n \leq 5 \cdot 2^t$. □

LEMMA 7. *For $N = 32 \cdot 2^t$, $t \geq 0$, and $N/4 + 1 \leq n \leq 9N/32$, projection designs of the maximal regular design of resolution IV with $9N/32$ factors do not have minimum aberration.*

PROOF. It is sufficient to show that the upper bound (with $N = 16 \cdot 2^{t+1}$) given in Corollary 2 is less than the lower bound given in Lemma 4, that is,

$$U(n) = [(65 \cdot 2^{3t+1} - 75 \cdot 2^{2t} + 5 \cdot 2^t)/6] \binom{n}{4} \bigg/ \binom{5 \cdot 2^{t+1}}{4} < L(n),$$

where $L(n)$ is defined in Lemma 4. Let $T = 2^t$. Given $t \geq 0$, $F(n) = L(n) - U(n)$ is a polynomial of degree 4. It is sufficient to show that $F(n)$ increases strictly when $8T \leq n \leq 9T$ and that $F(8T + 1) > 0$. Let $F^{(i)}(n)$ be the $i$th derivative of $F(n)$. It is straightforward to verify that $F^{(3)}(9T) > 0$ and $F^{(4)}(n) < 0$ when $8T \leq n \leq 9T$; therefore, $F^{(3)}(n) > 0$ for $8T \leq n \leq 9T$. It is easy to verify that $F^{(2)}(8T) > 0$, which implies $F^{(2)}(n) > 0$ for $8T \leq n \leq 9T$. It is also easy to verify that $F^{(1)}(8T) > 0$, leading to $F^{(1)}(n) > 0$ for $8T \leq n \leq 9T$. Therefore, $F(n)$ increases strictly when $8T \leq n \leq 9T$. Finally, it can be verified that $F(8T+1) > 0$. This completes the proof. □

THEOREM 3. *For $N = 32 \cdot 2^t$, $t \geq 0$, and $17N/64 \leq n \leq 5N/16$, a minimum aberration design with $N$ runs and $n$ factors must be a projection of the maximal regular design of resolution IV with $5N/16$ factors.*

PROOF. As explained in the Introduction, for $17N/64 \leq n \leq 5N/16$, every regular design must be a projection of maximal regular designs of resolution IV with $N/2$, $5N/16$, $9N/32$ or $17N/64$ factors. Lemmas 6 and 7 state that projection designs of maximal regular designs of resolution IV with $N/2$ or $9N/32$ factors do not have minimum aberration. We need only show that maximal designs of resolution IV with $17N/64$ factors do not have minimum aberration. According to Chen and Cheng [4], there are five maximal designs of resolution IV with $17N/64$ factors that are constructed by repeatedly doubling five maximal regular $64 \times 17$ designs of resolution IV. According to Block and Mee [2], all five maximal regular $64 \times 17$ designs do not have minimum aberration. Then it follows from Corollary 2.4 of [4] that all five maximal regular designs of resolution IV with $17N/64$ factors do not have minimum aberration. This completes the proof. □

Theorems 2 and 3 together lead to the following result on minimum aberration designs.

THEOREM 4. *Let $D$ and $\overline{D}$ be a pair of complementary projection designs of the maximal regular design of resolution IV with $5N/16$ factors. For $N = 32 \cdot 2^t$, $t \geq 0$, and $17N/64 \leq n \leq 5N/16$, $D$ has minimum aberration among all regular designs with $N$ runs and $n$ factors if $\overline{D}$ satisfies the following two conditions:*

(i) $A_4(\overline{D}), -A_5(\overline{D}), A_6(\overline{D}), -A_7(\overline{D}), \ldots$ *are sequentially minimized among all $u$-factor projection designs of the maximal regular design with $5N/16$ factors, where $u = 5N/16 - n$.*

(ii) $|f_i - f_j| \leq 1$ *for all $1 \leq i < j \leq 5$.*

As a numerical illustration, consider the maximal regular design of resolution IV with 256 runs and 80 factors. In Section 4 we constructed complementary designs $\overline{D}$ satisfying both conditions (i) and (ii) for $u = 1$–11. By Theorem 4, deleting them yields minimum aberration designs among all regular designs with $n = 69$–79 factors. Block [1] previously considered the construction of minimum aberration designs from this maximal regular design via a naive projection (i.e., deleting one column at a time). Our theoretical result confirms that the designs given by Block [1] with $n = 69$–$79$ factors have minimum aberration except for $n = 71$.

**Acknowledgments.** The research was initiated when the first author was visiting the Institute of Statistical Science, Academia Sinica. The authors thank an Associate Editor and a referee for their helpful comments.

## REFERENCES

[1] BLOCK, R. M. (2003). Theory and construction methods for large regular resolution IV designs. Ph.D. dissertation, Univ. Tennessee, Knoxville.
[2] BLOCK, R. M. and MEE, R. W. (2003). Second order saturated resolution IV designs. *J. Statist. Theory Appl.* **2** 96–112. MR2040435
[3] BUTLER, N. A. (2003). Some theory for constructing minimum aberration fractional factorial designs. *Biometrika* **90** 233–238. MR1966563
[4] CHEN, H. and CHENG, C.-S. (2006). Doubling and projection: A method of constructing two-level designs of resolution IV. *Ann. Statist.* **34** 546–558. MR2275253
[5] CHEN, H. and CHENG, C.-S. (2006). Some results on $2^{n-m}$ designs of resolution IV with (weak) minimum aberration. Preprint.
[6] CHEN, H. and HEDAYAT, A. S. (1996). $2^{n-l}$ designs with weak minimum aberration. *Ann. Statist.* **24** 2536–2548. MR1425966
[7] FRIES, A. and HUNTER, W. G. (1980). Minimum aberration $2^{k-p}$ designs. *Technometrics* **22** 601–608. MR0596803
[8] PLESS, V. (1963). Power moment identities on weight distributions in error correcting codes. *Inform. and Control* **6** 147–152. MR0262004
[9] RAKTOE, B. L., HEDAYAT, A. S. and FEDERER, W. T. (1981). *Factorial Designs*. Wiley, New York. MR0633756
[10] TANG, B. and WU, C. F. J. (1996). Characterization of minimum aberration $2^{n-m}$ designs in terms of their complementary designs. *Ann. Statist.* **24** 2549–2559. MR1425967
[11] WU, C. F. J. and HAMADA, M. (2000). *Experiments*: *Planning, Analysis and Parameter Design Optimization*. Wiley, New York. MR1780411
[12] XU, H. (2003). Minimum moment aberration for nonregular designs and supersaturated designs. *Statist. Sinica* **13** 691–708. MR1997169
[13] XU, H. (2005). A catalogue of three-level regular fractional factorial designs. *Metrika* **62** 259–281. MR2274993
[14] XU, H. (2006). Blocked regular fractional factorial designs with minimum aberration. *Ann. Statist.* **34** 2534–2553. MR2291509

DEPARTMENT OF STATISTICS
UNIVERSITY OF CALIFORNIA
LOS ANGELES, CALIFORNIA 90095-1554
USA
E-MAIL: hqxu@stat.ucla.edu

DEPARTMENT OF STATISTICS
UNIVERSITY OF CALIFORNIA
BERKELEY, CALIFORNIA 94720-3860
USA
E-MAIL: cheng@stat.berkeley.edu